\documentclass{amsart}
\usepackage{amsmath}
\usepackage{amssymb}
\usepackage{graphicx}
\usepackage{todonotes}
\usepackage{amsaddr}
\usepackage[a4paper, total={6in, 8in}]{geometry}
\usepackage{xcolor}
\usepackage{soul}
\usepackage{comment}

\newcommand{\ri}{\mathrm{i}}

\newtheorem{theorem}{Theorem}[section]
\theoremstyle{remark}
\newtheorem{remark}[theorem]{Remark}

\usepackage[colorlinks,linkcolor=red,citecolor=blue,urlcolor=blue]{hyperref}

\title{{A Note on Non-Hydrodynamic Solutions of Kinetic Systems}}

\author{Florian Kogelbauer}
\email{floriank@ethz.ch}
\address{Department of Mechanical and Process Engineering, ETH Z{u}rich, 8092 Z{u}rich, Switzerland}

\author{Ilya Karlin}
\email{ikarlin@ethz.ch}
\address{Department of Mechanical and Process Engineering, ETH Z{u}rich, 8092 Z{u}rich, Switzerland}

\date{\today}

\begin{document}

\begin{abstract}
    We show that the one-dimensional three-component Grad system admits solutions that violate the Chapman--Enskog scaling in Knudsen number. In particular, there exist solutions that do not converge to the analogues of the Euler and Navier--Stokes equations for vanishing Knudsen number. These non-hydrodynamic solutions correspond to a fast spectral manifold in kinetic phase space. 
\end{abstract}

 \maketitle

\section{Introduction}

Hilbert \cite{hilbert2022mathematical} famously posed the problem of connecting kinetic theory to continuum mechanics as the sixth of his list of mathematical challenges for the twentieth century. A modern interpretation of Hilbert's sixth problem entails showing that solutions of the Boltzmann equation converge - or at least stay close - to solutions to the Navier--Stokes equation in the limit of vanishing Knudsen number for some time \cite{saint2009hydrodynamic}. There is a huge body of work related to the convergence of scaled solutions to the Boltzmann equation, from which we can cite only  few \cite{bardos1991classical,caflisch1980fluid,de1989incompressible,gallagher2020convergence,guo2006boltzmann,nishida1978fluid}. Recently, the hard-sphere Boltzmann equation has been derived from microscopic Newtonian dynamics for larger times \cite{deng2024long}, which has been used to show that local Maxwellians evaluated at solutions to the Navier--Stokes equations stay close to true solutions of the Boltzmann equation for some time \cite{deng2025hilbertssixthproblemderivation}.\\

The main motivation for research along scaled solutions and convergence towards the Navier--Stokes equations is derived from the Chapman--Enskog series \cite{chapman1990mathematical} - a power series expansion in the Knudsen number. Indeed, at zeroth order, the Chapman--Enskog expansion formally recovers the Euler equations, and at first order the Navier--Stokes equations of fluid dynamics. Most of the research around Hilbert's Sixth Problem was thus focused on showing a relation between the Boltzmann and the Navier--Stokes and Euler equations. There are, however, several issues with the Chapman--Enskog expansion as a power-series expansion for a singular perturbation parameter.\\

Indeed, higher-order contributions of the Chapman--Enskog series can lead to non-physical behavior, such as the infamous Bobylev instability \cite{bobylev2006instabilities} for the Burnett equation. Furthermore, it is well-known that the limit from the Boltzmann equation to the Navier--Stokes dynamics is not unique as demonstrated, e.g., for the linear Boltzmann equation in \cite{ellis1974asymptotic}. It is thus questionable to assume that the Navier--Stokes dynamics will serve as a somewhat distinguished and universal fluid model derived from the Boltzmann equation. \\

In this work, we put forward a different interpretation of the relation between hydrodynamics and kinetic theory, based on the following hypothesis: Given a kinetic model for the time-evolution of a distribution function, such as the Boltzmann equation, hydrodynamics correspond to an invariant manifold in the phase space of distribution functions, called the \textit{hydrodynamic manifold}. This notion of hydrodynamics was already insinuated by Hilbert and is related to the concept of \textit{grossly determined solutions} \cite{truesdell1981fundamentals}. For more details and a broader overview in this research direction, we refer to \cite{Gorban01041994,gorban2005invariant,gorban2014hilbert}.  Recently, the notion of slow spectral manifold \cite{CABRE2005444} in conjunction with the spectral theory of kinetic operators \cite{ellis1975first} has been successfully applied to uniquely identify and explicitly describe the hydrodynamic manifold for linear kinetic systems \cite{kogelbauer2024exact,PhysRevE.110.055105}. \\

The ordering of eigenvalues of a kinetic system by negative real parts induces a hierarchy of time scales, which in turn corresponds to a nested family of slow invariant manifolds.  The remainder of the spectrum - usually the essential spectrum - corresponds to the fluctuation part of the kinetic dynamics and defines the fast, quickly decaying parts of solutions. The existence of the hydrodynamic manifold as a slow manifold leads to a distinguished, unique closure relation, i.e., a constitutive law that expresses higher-order fluxes in terms of basic hydrodynamic quantities, see \cite{kogelbauer2025exact} for an explicit calculation of the spectral closure for the linear BGK model. The spectrally closed constitutive law is dynamically optimal, holds for any level of rarefaction and can be used to make accurate predictions about rarefied flows properties, such as light-scattering spectra \cite{kogelbauer2025learning} or rarefaction effects 
in shear flows \cite{4gkn-7s3x}, solely in terms of macroscopic field equations.\\

In this note, we are not concerned with the hydrodynamic manifold, but much rather with the non-hydrodynamic part of kinetic dynamics. In particular, we give an illustration of non-hydrodynamic kinetic solutions corresponding to fast dynamics that violate the Chapman--Enskog scaling. That is, we show that there exist solutions to the kinetic model whose hydrodynamic moments do not converge to the analogues of the Euler or the Navier--Stokes dynamics for vanishing Knudsen number.  We focus on a linear, one-dimensional, three-component Grad system, which belongs to the widely used class of moment projections of kinetic models. For this model, the Chapman--Enskog series can be summed in closed form \cite{karlin2002hydrodynamics} and a viscosity-capillarity balance was derived \cite{slemrod2012chapman}. It was shown in \cite{kogelbauer2020slow} that the slow modes at each wave number span a two-dimensional hydrodynamic manifold. The solutions on this manifold converge to the analogues of the Euler and Navier--Stokes equations for vanishing Knudsen numbers. The non-hydrodynamic solutions then correspond to  the remaining fast manifold at each wave number, for which we show that they diverge in Knudsen number. This divergence in decay rates gives further indication that the dynamic, short-time proximity of certain solutions to the Boltzmann equation to solutions of the Navier--Stokes equations might not hold true for large times or a larger subset of distribution functions.

\section{The Three-Component Grad System: Spectral Analysis in Frequency Space}

\begin{figure}
    \centering
\includegraphics[width=0.5\linewidth]{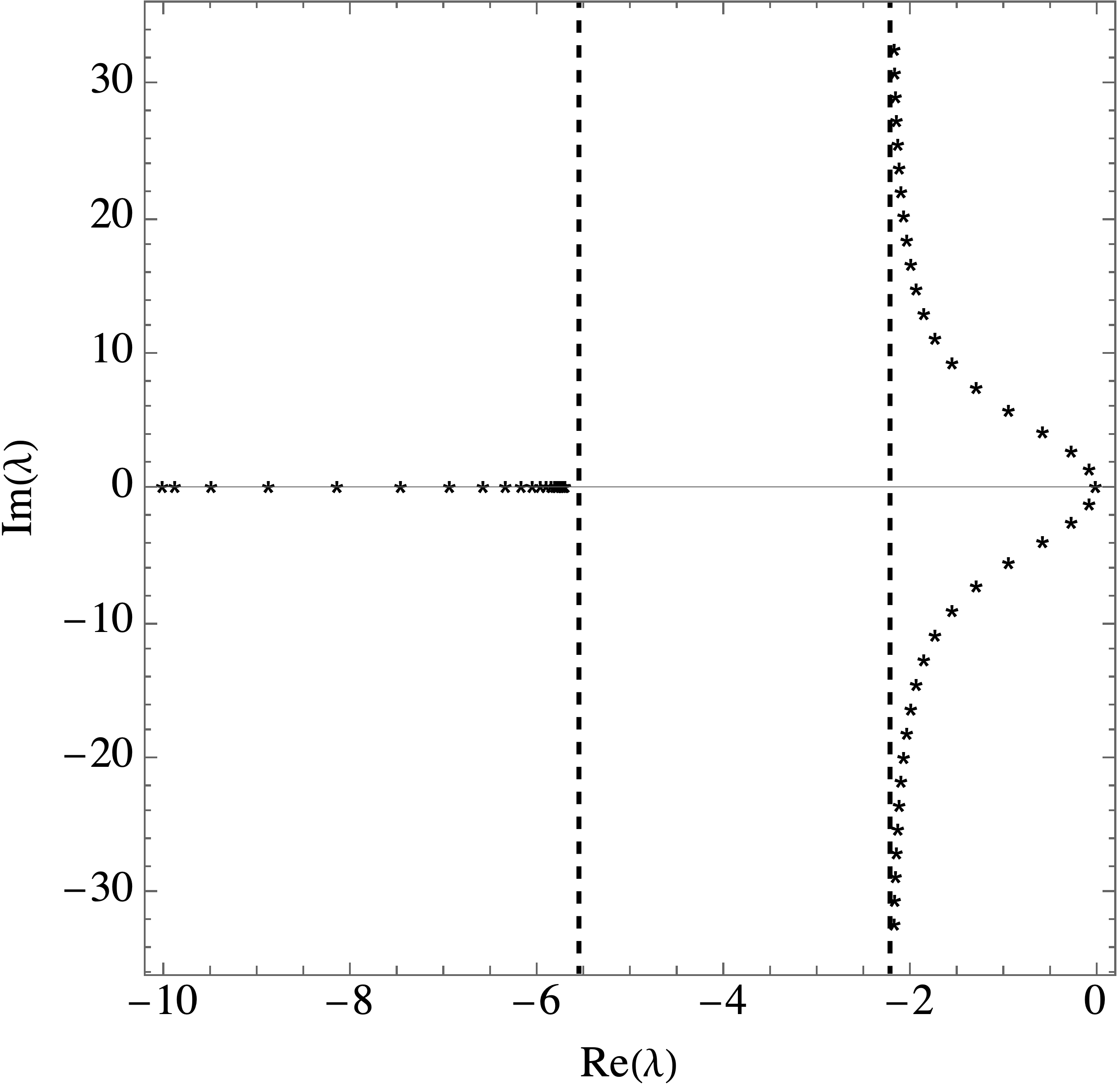}
    \caption{
    Spectrum of the  linear Grad system \eqref{eqGrad} for $\varepsilon=0.1$ and various wave numbers as derived in, e.g., \cite{kogelbauer2020slow}. The fast real modes (star symbol on the real axis) accumulate at $-\frac{5}{9\varepsilon}$ (vertical dashed line at $\Re{\lambda} = -50/9$), while the slow complex conjugated modes (bell-shaped star symbol) accumulate at $\Re\lambda = -\frac{2}{9\varepsilon}$ (dashed line at $\Re{\lambda} = -20/9$) as $k\to\infty$. The smaller the Knudsen number, the more the axis of accumulation of the fast modes moves the left. }
    \label{pspec_Grad3D}
\end{figure}

Consider the three-component Grad system,
\begin{equation}\label{eqGrad}
\begin{split}
&\frac{\partial p}{\partial t}=-\frac{5}{3}\frac{\partial u}{\partial x},\\
&\frac{\partial u}{\partial t}=-\frac{\partial p}{\partial x}-\frac{\partial \sigma}{\partial x},\\
&\frac{\partial \sigma}{\partial t}=-\frac{4}{3}\frac{\partial u}{\partial x}-\frac{1}{\varepsilon}\sigma,
\end{split}
\end{equation}
for the pressure $p$, the velocity $u$ and the stress $\sigma$. The parameter $\epsilon>0$ represents the Knudsen number, i.e., the ratio of the mean-free path to a typical length scale of the system. Equation \eqref{eqGrad} is derived from a 13-moment system by setting the heat flux to zero \cite{grad1949kinetic} and serves as a toy model for kinetic dynamics. While \eqref{eqGrad} retains important properties of realistic kinetic models, such as global decay towards equilibrium due to the existence of an entropy, the Grad system under consideration does not exhibit criticality as the acoustic eigenvalues exist for all wave numbers.\\
In the following, we pose equation \eqref{eqGrad} on the whole real line with suitable decay at infinity. Periodic boundary conditions can be treated in a similar fashion. We denote the Fourier transform of a square-integrable function $f:\mathbb{R}\to\mathbb{R}$ as 
\begin{equation}
    \hat{f}(k) = \frac{1}{\sqrt{2\pi}} \int_{-\infty}^{\infty} f(x) e^{-\ri xk}\, dx .
    \end{equation}
In frequency space, system \eqref{eqGrad} reads,
\begin{equation}\label{eqGradk}
\frac{d}{d t}\left(\begin{array}{c}
     \hat{p}\\
     \hat{u}\\
     \hat{\sigma}
\end{array}\right) = L_k\left(\begin{array}{c}
     \hat{p}\\
     \hat{u}\\
     \hat{\sigma}
\end{array}\right),
\end{equation}
for the wave-number dependent matrix,
\begin{equation}\label{Lk}
    L_k = \begin{pmatrix}
        0 &  -\frac{5}{3}\ri k & 0\\
-\ri k & 0 & -\ri k\\
0 & -\frac{4}{3}\ri k & -\frac{1}{\varepsilon}
    \end{pmatrix},\quad k\in\mathbb{R}. 
\end{equation}
Let us recall the spectral properties of $L_k$ as derived in \cite{kogelbauer2020slow,slemrod2012chapman}. The characteristic polynomial of $L_k$
\begin{equation}\label{charAk}
P_{k}(\lambda)=-\lambda^3-\frac{1}{\varepsilon}\lambda^2-3k^2\lambda-\frac{5}{3}\frac{k^2}{\varepsilon},
\end{equation}
has a pair of complex conjugated eigenvalues $\{\lambda_{\rm  ac},\lambda_{\rm  ac}^*\}$, called acoustic modes, and one real eigenvalue $\lambda_{\rm diff}$, called generalized (non-hydrodynamic) diffusion mode.
The eigenvectors of $L_k$ are given by
\begin{equation}\label{defQ}
Q = \left(\begin{matrix}
-1-a_1b_1 & -1-a_2b_2 & -1-a_3b_3\\
\ri b_1 & \ri b_2 & \ri b_3\\
1 & 1 & 1
\end{matrix}\right),
\end{equation}
for the constants 
\begin{equation}
a_j=\frac{\lambda_j}{k},\ b_j=\frac{3}{4\varepsilon k}(1+\varepsilon \lambda_j),\ j=1,2,3, 
\end{equation}
where $\lambda_1 = \lambda_{\rm ac}$, $\lambda_2 = \lambda_{\rm ac}^*$ and $\lambda_{3} = \lambda_{\rm diff}$. A typical spectrum of $L_k$ ranging over different values of wave numbers is depicted in Figure \ref{pspec_Grad3D}. We emphasize that the eigenvectors in \eqref{defQ} are not orthogonal to each other since the matrix $L_k$ is not normal.

\begin{figure}
    \centering
\includegraphics[width=0.7\linewidth]{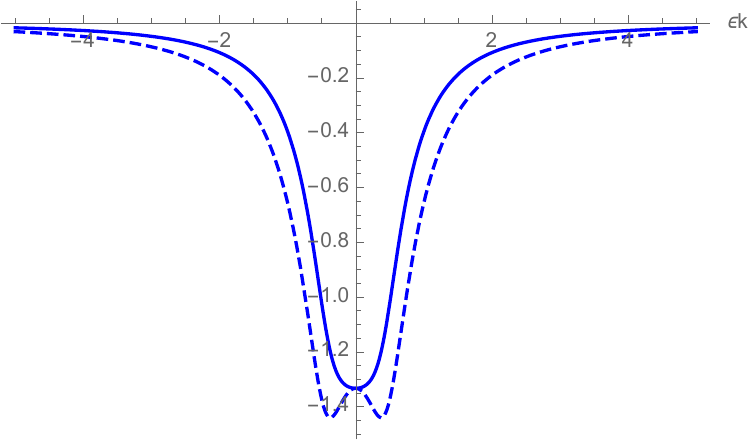}
    \caption{The functions $A$ (solid) and $B$ (dashed) in dependence of $\varepsilon k$. }
    \label{ABFigure}
\end{figure}

The slow manifold at each wave number is then given explicitly as the plane spanned by the two acoustic modes, i.e., the first two columns of \eqref{defQ},
\begin{equation}\label{Mslow}
    \mathcal{M}_{\rm slow} = \text{span} \left\{\left(\begin{array}{c}
         -1-a_1b_1 \\
         \ri b_1 \\
         1\end{array}\right),\left(\begin{array}{c}
         -1-a_2b_2 \\
         \ri b_2 \\
         1
    \end{array}\right)\right\}.
\end{equation}
On the slow manifold, we can close the $(p,u)$-dynamics uniquely by the relation
\begin{equation}\label{sigma}
\hat{\sigma} = \ri kA(k^2,\varepsilon)\hat{u} -k^2B(k^2,\varepsilon)\hat{p} ,
\end{equation}
where the real-valued functions $(k,\varepsilon)\mapsto A(k^2,\varepsilon)$ and $(k,\varepsilon)\mapsto B(k^2,\varepsilon)$ take the form
\begin{equation}\label{defAB}
\begin{split}
&A(k^2,\varepsilon)=-\frac{4\varepsilon(1+\varepsilon(\lambda_1+\lambda_2))}{3+3\varepsilon(\lambda_1+\lambda_2)+\varepsilon
^2(3\lambda_1\lambda_2-4k^2)},\\
&B(k^2,\varepsilon)=-\frac{4\varepsilon^2}{3+3\varepsilon(\lambda_1+\lambda_2)+\varepsilon
	^2(3\lambda_1\lambda_2-4k^2)}.
\end{split}
\end{equation}
The functions $A$ and $B$ are depicted in Figure \ref{ABFigure}.\\

We call the one-dimensional linear sub-space spanned by the remaining eigenvector \textit{non-hydrodynamic manifold} or \textit{fast manifold}:
\begin{equation}\label{Mfast}
 \mathcal{M}_{\rm fast} = \text{span} \left\{\left(\begin{array}{c}
         -1-a_3b_3 \\
         \ri b_3 \\
         1\end{array}\right)
         \right\}. 
\end{equation}

Analogous to the slow manifold closure \eqref{sigma}, the fast manifold induces a constitutive law through the relations
\begin{equation}\label{sigmafast}
    \hat{\sigma} = \frac{1}{\ri b_3}\hat{u} = -\frac{1}{1+a_3b_3}\hat{p}. 
\end{equation}

\begin{remark}
    The exact hydrodynamic closure for the Grad system was obtained in \cite{gorban1992structure} by full summation of the Chapman--Enskog series. In \cite{kogelbauer2020slow}, it was shown that the exact summation of the Chapman--Enskog series for the three-component Grad system is equivalent to the spectral closure obtained by projection onto the two-dimensional slow manifold \eqref{Mslow}. 
\end{remark}

\begin{remark}
    The idea that a Taylor series expansion in Knudsen number selects a special analytic invariant manifold which serves as a hydrodynamic manifold dates at least back to \cite{mckean1969simple}. Apart from the Bobylev instability for higher-order fluid models in Knudsen number, however, there might be even more severe problems with this approach. Indeed, even for simple models, the Chapman--Enskog series can be divergent altogether \cite{kogelbauer2025relationexacthydrodynamicschapmanenskog} and it is thus not clear if the summation of the Chapman--Enskog series always corresponds to the hydrodynamic manifold. The identification of the hydrodynamic manifold as a slow spectral manifold circumvents this issue by selecting the unique, invariant manifold corresponding to eigenvalues with minimal negative real part. 
\end{remark}

\section{Hydrodynamics on the Slow Manifold and Divergence in Decay Rates}
In this section, we first recall the construction of the spectral closure on the hydrodynamic manifold and show that the dynamics for $u$ and $p$ converge to the analogues of the Euler and Navier--Stokes equations in the limit of vanishing Knudsen number. We then show that the closure on the fast eigenmode leads to a different scaling in $\varepsilon$.\\

Plugging \eqref{sigma} into \eqref{eqGradk} leads to the hydrodynamics on the slow manifold:
\begin{equation}\label{slowhydro}
    \frac{\partial}{\partial t} \left(\begin{array}{c}
    \hat{p}_{s} \\ \hat{u}_{s}
    \end{array}\right) = \begin{pmatrix}
        0 & -\frac{5}{3}\ri k \\
        -\ri k (1-k^2B)  & k^2A
    \end{pmatrix} \left(\begin{array}{c}
    \hat{p}_{s} \\ \hat{u}_s
    \end{array}\right).
\end{equation}
Using the cubic solution formula, we find that the diffusion mode takes the explicit form,
\begin{equation}\label{lambda3explicit}
\begin{split}
	\lambda_3&=\frac{1}{3\varepsilon}\left[\sqrt[3]{-1-9k^2\varepsilon^2+3\varepsilon\sqrt{5k^2-18k^4\varepsilon^2+81k^6\varepsilon^4}}\right.\\
	&\left.\qquad+\sqrt[3]{-1-9k^2\varepsilon^2-3\varepsilon\sqrt{5k^2-18k^4\varepsilon^2+81k^6\varepsilon^4}}-1\right],
\end{split}
\end{equation}
which implies that 
\begin{equation}
    \lambda_{3}(k) = -\frac{1}{\varepsilon} + \mathcal{O}(1),\quad \varepsilon \to 0. 
\end{equation}
Consequently, we have that for fixed wave number,
\begin{equation}
\begin{split}
       A & = \varepsilon + \mathcal{O}(\varepsilon^2), \quad \varepsilon \to 0 \\
       B & = \varepsilon^2 + \mathcal{O}(\varepsilon^3), \quad \varepsilon \to 0, 
\end{split}
\end{equation}
and the Chapman--Enskog expansion of the stress,
\begin{equation}
    \hat{\sigma} = -\frac{4}{3}\varepsilon\ri k \hat{u} +\frac{4}{3}\varepsilon^2k^2\hat{p} + \mathcal{O}(\varepsilon^3),
\end{equation}
which recovers the analogue of the Euler equations on the slow manifold:
\begin{equation}
    \frac{\partial}{\partial t} \left(\begin{array}{c}
    \hat{p}_{s} \\ \hat{u}_{s}
    \end{array}\right) = \begin{pmatrix}
        0 & -\frac{5}{3}\ri k \\
        -\ri k   & 0
    \end{pmatrix} \left(\begin{array}{c}
    \hat{p}_{s} \\ \hat{u}_s
    \end{array}\right) + \mathcal{O}(\varepsilon). 
\end{equation}
At order $\varepsilon$, we find the analogue of the Navier--Stokes equations,
\begin{equation}
    \frac{\partial}{\partial t} \left(\begin{array}{c}
    \hat{p}_{s} \\ \hat{u}_{s}
    \end{array}\right) = \begin{pmatrix}
        0 & -\frac{5}{3}\ri k \\
        -\ri k   & -\frac{4}{3}\varepsilon
    \end{pmatrix} \left(\begin{array}{c}
    \hat{p}_{s} \\ \hat{u}_s
    \end{array}\right) + \mathcal{O}(\varepsilon^2).
\end{equation}
The previous considerations about the exact non-local hydrodynamics and their consistency to the Euler and Navier--Stokes equations through the Chapman--Enskog expansion are well established \cite{gorban2005invariant,kogelbauer2020slow}. We now show that the Grad system admits solutions that do not converge to a classical hydrodynamics, i.e., diverge as $\varepsilon\to 0$. This is achieved by considering the hydrodynamic closure on the fast manifold \eqref{Mfast} using the constitutive law \eqref{sigmafast}. On the one hand, those solutions are somewhat nonphysical, since a generic initial condition will converge exponentially fast to the two-dimensional slow manifold. On the other hand, the existence of an invariant manifold on which hydrodynamics diverge shows that the Chapman--Enskog scaling does not cover the full phase space.\\
Analogously, we obtain the hydrodynamics on the fast manifold spanned by the diffusion mode as 
\begin{equation}\label{fasthydro1}
    \frac{\partial}{\partial t}  \left(\begin{array}{c}
    \hat{p}_{f} \\ \hat{u}_{f}
    \end{array}\right) = \lambda_{3} \left(\begin{array}{c}
    \hat{p}_{f} \\ \hat{u}_{f}
    \end{array}\right),
\end{equation}
through the constitutive law \ref{sigmafast}.\\
Using the relation 
\begin{equation}\label{fastrelation}
    \lambda_{3} \hat{\rho}_f = -\frac{5}{3}\ri k \hat{u}_f,
\end{equation}
which follows from the eigenstructure \eqref{defQ} together with characteristic polynomial \eqref{charAk}, we can rewrite \eqref{fasthydro1} as 
\begin{equation}\label{fasthydro2}
        \frac{\partial}{\partial t}  \left(\begin{array}{c}
    \hat{p}_{f} \\ \hat{u}_{f}
    \end{array}\right) = \begin{pmatrix}
       0 &  -\frac{5}{3}\ri k \\
       0 & \lambda_3 
    \end{pmatrix} \left(\begin{array}{c}
    \hat{p}_{f} \\ \hat{u}_{f}
    \end{array}\right).
\end{equation}
We now claim that the fast hydrodynamics \eqref{fasthydro2} are divergent in the limit
as $\varepsilon\to 0$ point-wise in time along their trajectories. To show this, we have to make sure that there is no linear combination of $\hat{p}_f$ and $\hat{u}_f$ in the right-hand side of \eqref{fasthydro1} that remains bounded as $\varepsilon\to 0$. This is a consequence of the constitutive law \eqref{sigmafast}: The fast manifold is one-dimensional, while the hydrodynamics are written for the two fields $(\hat{p},\hat{u})$. Indeed, the $\hat{p}_f$-dynamics in \eqref{fasthydro1} are divergent for $\varepsilon\to 0$ when written in the $\hat{p}_f$-variable, but are bounded when written in the $\hat{u}_f$-variable.\\
We cannot rewrite \eqref{fasthydro2} using \eqref{fastrelation} to match the slow hydrodynamics \eqref{slowhydro}. We can expand 
\begin{equation}\label{fastexpand}
\begin{split}
      \lambda_3 \hat{u}_f & = [\alpha + (\lambda_3-\alpha)]\hat{u}_f\\
      &  = \frac{3\ri}{5k}\alpha \lambda_3\hat{p}_f + (\lambda_3-\alpha)\hat{u}_f,
\end{split}
\end{equation}
for some function $\alpha$. To make the second term on the right-hand side in \eqref{fastexpand} vanish, we have to have that $\alpha = -1/\varepsilon+\mathcal{O}(1)$ (with actually the order-one terms matching as well), while for the first term to be of order one, we have to have that $\alpha = \mathcal{O}(\varepsilon)$ - a contradiction. This shows that any hydrodynamics on the fast manifold will diverge in the limit of vanishing Knudsen number.\\


\begin{remark}
    Slemrod \cite{slemrod2012chapman} demonstrated  that the the exact summation of the Chapman--Enskog series for the three-component Grad system implies a viscosity-capillarity balance of the form
    \begin{equation}\label{visccap}
        \frac{1}{2}\frac{\partial}{\partial t} \int_{-\infty}^{\infty} \frac{3}{5}|p|^2+|u|^2\, dx + \frac{\partial}{\partial t} \int_{-\infty}^\infty -\frac{5}{3}k^2 B(k^2) |\hat{p}|^2\, dk = \int_{-\infty}^\infty k^2 A(k^2) |\hat{u}|^2\,dk. 
    \end{equation}
The existence of such a balance law was extended to general kinetic models in \cite{PhysRevE.110.055105}. Relation \eqref{visccap} can be interpreted as a non-local generalization of Korteweg's capillary action. Slemrod argued \cite{slemrod2013boltzmann} that the existence of balance law of the form \eqref{visccap} might pose an obstruction to the global convergence of kinetic solutions to solutions of the Navier--Stokes equations in the nonlinear regime. The intuition is based on an analogy to the convergence properties of the Korteweg–de Vries–Burgers equations for vanishing viscosity and capillarity. This would imply a different limiting behavior on the (nonlinear) hydrodynamic manifold. In this note, we have shown the somewhat more modest divergence on the fast manifold, but even in the linear regime. 
\end{remark}

\begin{remark}
We mention that the special fast solutions of the Grad system as constructed in \eqref{fasthydro2} cannot be identified through a perturbative solution of the invariance equation directly. Indeed, iterative schemes typically solve for the slow, hydrodynamic manifold due to implicit numerical stability \cite{colangeli2007hyperbolicity}. The fast non-hydrodynamic manifold can thus usually not be detected through a generic numerical investigation of the invariance equation and a detailed spectral analysis is necessary to infer a divergence in decay rates for non-hydrodynamic solutions. 
\end{remark}

\section{Conclusions}

In this note, we demonstrated that the one-dimensional three-component Grad system admits non-hydrodynamic solutions whose behavior fundamentally deviates from the predictions of the Chapman–Enskog expansion in the limit of vanishing Knudsen number. By identifying the fast manifold associated with the non-hydrodynamic diffusion mode, we provided a clear spectral characterization of solutions that do not converge to the Euler or Navier–Stokes equations. Our analysis highlights the limitations of traditional hydrodynamic closures based on Chapman–Enskog scaling and supports an alternative view of hydrodynamics as dynamics confined to a slow invariant manifold. These findings suggest that, even in simple kinetic models, the distinction between fast and slow modes is not merely technical but has profound implications for the structure and universality of macroscopic fluid dynamics derived from kinetic theory. It is not clear at this point whether the existence of non-hydrodynamic modes that violate the Chapman--Enskog scaling are a feature of the moment-truncation employed in deriving the Grad system, or if it is an inherent feature of kinetic models. It would thus be interesting to investigate the possibility of non-hydrodynamic solutions for more realistic kinetic models, such as the one-dimensional linear kinetic model for which a complete and explicit spectral theory is known \cite{kogelbauer2021non}.

\section{Acknowledgment}

The authors would like to thank Marshall Slemrod for several useful comments on the manuscript.

\bibliographystyle{abbrv}
\bibliography{DivergenceinDecay Rates}

\end{document}